# Heisenberg-Pauli-Weyl and Donoho Stark's uncertainty principle for the Weinstein $L^2$-multiplier operators

Ahmed Saoudi

**Abstract.** The aim of this paper is establish the Heisenberg-Pauli-Weyl uncertainty principle and Donho-Stark's uncertainty principle for the Weinstein $L^2$-multiplier operators.

## 1. Introduction

The Weinstein operator $\Delta_{W,\alpha}^d$ defined on $\mathbb{R}_+^{d+1} = \mathbb{R}^d \times (0, \infty)$, by

$$\Delta_{W,\alpha}^d = \sum_{j=1}^{d+1} \frac{\partial^2}{\partial x_j^2} + \frac{2\alpha+1}{x_{d+1}} \frac{\partial}{\partial x_{d+1}} = \Delta_d + L_\alpha, \ \alpha > -1/2,$$

where $\Delta_d$ is the Laplacian operator for the $d$ first variables and $L_\alpha$ is the Bessel operator for the last variable defined on $(0, \infty)$ by

$$L_\alpha u = \frac{\partial^2 u}{\partial x_{d+1}^2} + \frac{2\alpha+1}{x_{d+1}} \frac{\partial u}{\partial x_{d+1}}.$$

The Weinstein operator $\Delta_{W,\alpha}^d$ has several applications in pure and applied mathematics, especially in fluid mechanics [4].

The Weinstein transform generalizing the usual Fourier transform, is given for $\varphi \in L_\alpha^1(\mathbb{R}_+^{d+1})$ and $\lambda \in \mathbb{R}_+^{d+1}$, by

$$\mathcal{F}_{W,\alpha}(\varphi)(\lambda) = \int_{\mathbb{R}_+^{d+1}} \varphi(x) \Lambda_\alpha^d(x, \lambda) d\mu_\alpha(x),$$

where $d\mu_\alpha(x)$ is the measure on $\mathbb{R}_+^{d+1} = \mathbb{R}^d \times (0, \infty)$ and $\Lambda_\alpha^d$ is the Weinstein kernel given respectively later by (2.1) and (2.4).

Let $m$ be a function in $L^2_\alpha(\mathbb{R}^{d+1}_+)$ and let $\sigma$ be a positive real number. The Weinstein $L^2$-multiplier operators is defined for smooth functions $\varphi$ on $\mathbb{R}^{d+1}_+$, in [14] as

$$\mathcal{T}_{w,m,\sigma}\varphi(x) := \mathcal{F}^{-1}_{W,\alpha}\left(m_\sigma \mathcal{F}_{W,\alpha}(\varphi)\right)(x), \quad x \in \mathbb{R}^{d+1}_+, \tag{1.1}$$

where the function $m_\sigma$ is given by

$$m_\sigma(x) = m(\sigma x).$$

These operators are a generalization of the multiplier operators $\mathcal{T}_m$ associated with a bounded function $m$ and given by $\mathcal{T}_m(\varphi) = \mathcal{F}^{-1}(m\mathcal{F}(\varphi))$, where $\mathcal{F}(\varphi)$ denotes the ordinary Fourier transform on $\mathbb{R}^n$. These operators gained the interest of several Mathematicians and they were generalized in many settings in [1, 3, 6, 13, 14, 16, 17, 18].

In this work we are interested the $L^2$ uncertainty principles for the Weinstein multiplier operators. The uncertainty principles play an important role in harmonic analysis. These principles state that a function $\varphi$ and its Fourier transform $\mathcal{F}(\varphi)$ cannot be simultaneously sharply localized. Many aspects of such principles are studied for several Fourier transforms.

Many uncertainty principles have already been proved for the Weinstein transform $\mathcal{F}_{W,\alpha}$, namely by N. Ben Salem, A. R. Nasr [2] and Mejjaoli H. and Salhi M. [9]. The authors have established in [9] the Heisenberg-Pauli-Weyl inequality for the Weinstein transform, by showing that, for every $\varphi$ in $L^2_\alpha(\mathbb{R}^{d+1}_+)$

$$\|\varphi\|_{\alpha,2} \le \frac{2}{2\alpha + d + 2} \||x|\varphi\|_{\alpha,2} \||y|\mathcal{F}_{W,\alpha}(\varphi)\|_{\alpha,2}. \tag{1.2}$$

In the present paper we are interested in proving an analogue of Heisenberg-Pauli-Weyl uncertainty principle For the operators $\mathcal{T}_{w,m,\sigma}$. More precisely, we will show, for $\varphi \in L^2_\alpha(\mathbb{R}^{d+1}_+)$

$$\|\varphi\|_{\alpha,2} \le \frac{2\||y|\mathcal{F}_{W,\alpha}(\varphi)\|_{\alpha,2}}{2\alpha + d + 2} \left( \int_{\mathbb{R}^{d+1}_+} \int_0^\infty |x|^2 |\mathcal{T}_{w,m,\sigma}\varphi(x)|^2_{\alpha,2} \frac{d\sigma}{\sigma} d\mu_\alpha(x) \right)^{\frac{1}{2}},$$

provided $m$ be a function in $L^2_\alpha(\mathbb{R}^{d+1}_+)$ satisfying the admissibility condition

$$\int_0^\infty |m_\sigma(x)| \frac{d\sigma}{\sigma} = 1, \quad \text{a.e. } x \in \mathbb{R}^{d+1}_+. \tag{1.3}$$

Moreover, for $\beta, \delta \in [1, \infty)$ and $\varepsilon \in \mathbb{R}$, such that $\beta\varepsilon = (1-\varepsilon)\delta$, we will show

$$\|\varphi\|_{\alpha,2} \le \left( \frac{2}{2\alpha + d + 2} \right)^{\beta\varepsilon} \||x|^\beta \mathcal{T}_{w,m,\sigma}\varphi\|^\varepsilon_{\alpha,2} \||y|^\delta \mathcal{F}_{W,\alpha}(\varphi)\|^{1-\varepsilon}_{\alpha,2}.$$

Using the techniques of Donoho and Stark [5], we show uncertainty principle of concentration type for the $L^2$ theory. Let $\varphi$ be a function in $L^2_\alpha(\mathbb{R}^{d+1}_+)$ and $m \in L^1_\alpha(\mathbb{R}^{d+1}_+) \cap L^2_\alpha(\mathbb{R}^{d+1}_+)$ satisfying the admissibility condition (1.3). If $\varphi$ is $\epsilon$-concentrated on $\Omega$ and $\mathcal{T}_{w,m,\sigma}\varphi$ is $\nu$-concentrated on $\Sigma$, then

$$\|m\|_{\alpha,1} (\mu_\alpha(\Omega))^{\frac{1}{2}} \left( \int \int_\Sigma \frac{1}{\sigma^{2(2\alpha+d+2)}} d\Theta_\alpha(\sigma, x) \right)^{\frac{1}{2}} \ge 1 - (\epsilon + \nu),$$

where $\Theta_\alpha$ is the measure on $(0,\infty) \times \mathbb{R}^{d+1}_+$ given by $d\Theta_\alpha(\sigma, x) := (d\sigma/\sigma) d_\alpha \mu(x)$.

This paper is organized as follows. In section 2, we recall some basic harmonic analysis results related with the Weinstein operator $\Delta^d_{W,\alpha}$ and we introduce preliminary facts that will be used later.

In section 3, we establish Heisenberg-Pauli-Weyl uncertainty principle For the operators $\mathcal{T}_{w,m,\sigma}$.

The last section of this paper is devoted to Donoho-Stark's uncertainty principle for the Weinstein $L^2$-multiplier operators.

## 2. Harmonic analysis Associated with the Weinstein Operator

In this section, we shall collect some results and definitions from the theory of the harmonic analysis associated with the Weinstein operator $\Delta^d_{W,\alpha}$. Main references are [10, 11, 12].

In the following we denote by

- $\mathbb{R}^{d+1}_+ = \mathbb{R}^d \times (0, \infty)$.
- $x = (x_1, ..., x_d, x_{d+1}) = (x', x_{d+1})$.
- $-x = (-x', x_{d+1})$.
- $C_*(\mathbb{R}^{d+1})$, the space of continuous functions on $\mathbb{R}^{d+1}$, even with respect to the last variable.
- $S_*(\mathbb{R}^{d+1})$, the space of the $C^\infty$ functions, even with respect to the last variable, and rapidly decreasing together with their derivatives.
- $L^p_\alpha(\mathbb{R}^{d+1}_+)$, $1 \leq p \leq \infty$, the space of measurable functions $f$ on $\mathbb{R}^{d+1}_+$ such that

$$\|f\|_{\alpha,p} = \left(\int_{\mathbb{R}^{d+1}_+} |f(x)|^p \, d\mu_\alpha(x)\right)^{1/p} < \infty, \ p \in [1, \infty),$$

$$\|f\|_{\alpha,\infty} = \operatorname*{ess\,sup}_{x \in \mathbb{R}^{d+1}_+} |f(x)| < \infty,$$

where

$$d\mu_\alpha(x) = \frac{x_{d+1}^{2\alpha+1}}{(2\pi)^d 2^{2\alpha} \Gamma^2(\alpha+1)} dx. \tag{2.1}$$

- $\mathcal{A}_\alpha(\mathbb{R}^{d+1}) = \{\varphi \in L^1_\alpha(\mathbb{R}^{d+1}_+); \mathcal{F}_{W,\alpha}\varphi \in L^1_\alpha(\mathbb{R}^{d+1}_+)\}$ the Wiener algebra space.

We consider the Weinstein operator $\Delta^d_{W,\alpha}$ defined on $\mathbb{R}^{d+1}_+$ by

$$\Delta^d_{W,\alpha} = \sum_{j=1}^{d+1} \frac{\partial^2}{\partial x_j^2} + \frac{2\alpha+1}{x_{d+1}} \frac{\partial}{\partial x_{d+1}} = \Delta_d + L_\alpha, \ \alpha > -1/2, \tag{2.2}$$

where $\Delta_d$ is the Laplacian operator for the $d$ first variables and $L_\alpha$ is the Bessel operator for the last variable defined on $(0, \infty)$ by

$$L_\alpha u = \frac{\partial^2 u}{\partial x_{d+1}^2} + \frac{2\alpha+1}{x_{d+1}} \frac{\partial u}{\partial x_{d+1}}.$$

The Weinstein operator $\Delta_{W,\alpha}^d$ have remarkable applications in diffrent branches of mathematics. For instance, they play a role in Fluid Mechanics [4].

### 2.1. The eigenfunction of the Weinstein operator

For all $\lambda = (\lambda_1, ..., \lambda_{d+1}) \in \mathbb{C}^{d+1}$, the system

$$\frac{\partial^2 u}{\partial x_j^2}(x) = -\lambda_j^2 u(x), \quad \text{if } 1 \leq j \leq d$$
$$L_\alpha u(x) = -\lambda_{d+1}^2 u(x), \quad (2.3)$$
$$u(0) = 1, \quad \frac{\partial u}{\partial x_{d+1}}(0) = 0, \quad \frac{\partial u}{\partial x_j}(0) = -i\lambda_j, \quad \text{if } 1 \leq j \leq d$$

has a unique solution denoted by $\Lambda_\alpha^d(\lambda, .)$, and given by

$$\Lambda_\alpha^d(\lambda, x) = e^{-i<x',\lambda'>} j_\alpha(x_{d+1}\lambda_{d+1}) \quad (2.4)$$

where $x = (x', x_{d+1})$, $\lambda = (\lambda', \lambda_{d+1})$ and $j_\alpha$ is is the normalized Bessel function of index $\alpha$ defined by

$$j_\alpha(x) = \Gamma(\alpha + 1) \sum_{k=0}^\infty \frac{(-1)^k x^{2k}}{2^k k! \Gamma(\alpha + k + 1)}.$$

The function $(\lambda, x) \mapsto \Lambda_\alpha^d(\lambda, x)$ has a unique extension to $\mathbb{C}^{d+1} \times \mathbb{C}^{d+1}$, and satisfied the following properties.

**Proposition 2.1.** *i). For all $(\lambda, x) \in \mathbb{C}^{d+1} \times \mathbb{C}^{d+1}$ we have*

$$\Lambda_\alpha^d(\lambda, x) = \Lambda_\alpha^d(x, \lambda). \quad (2.5)$$

*ii). For all $(\lambda, x) \in \mathbb{C}^{d+1} \times \mathbb{C}^{d+1}$ we have*

$$\Lambda_\alpha^d(\lambda, -x) = \Lambda_\alpha^d(-\lambda, x). \quad (2.6)$$

*iii). For all $(\lambda, x) \in \mathbb{C}^{d+1} \times \mathbb{C}^{d+1}$ we get*

$$\Lambda_\alpha^d(\lambda, 0) = 1. \quad (2.7)$$

*vi). For all $\nu \in \mathbb{N}^{d+1}$, $x \in \mathbb{R}^{d+1}$ and $\lambda \in \mathbb{C}^{d+1}$ we have*

$$\left| D_\lambda^\nu \Lambda_\alpha^d(\lambda, x) \right| \leq \|x\|^{|\nu|} e^{\|x\| \|\Im \lambda\|} \quad (2.8)$$

*where $D_\lambda^\nu = \partial^\nu / (\partial \lambda_1^{\nu_1} ... \partial \lambda_{d+1}^{\nu_{d+1}})$ and $|\nu| = \nu_1 + ... + \nu_{d+1}$. In particular, for all $(\lambda, x) \in \mathbb{R}^{d+1} \times \mathbb{R}^{d+1}$, we have*

$$\left| \Lambda_\alpha^d(\lambda, x) \right| \leq 1. \quad (2.9)$$

## 2.2. The Weinstein transform

**Definition 2.2.** The Weinstein transform is given for $\varphi \in L^1_\alpha(\mathbb{R}^{d+1}_+)$ by

$$\mathcal{F}_{W,\alpha}(\varphi)(\lambda) = \int_{\mathbb{R}^{d+1}_+} \varphi(x)\Lambda^d_\alpha(\lambda, x)d\mu_\alpha(x), \quad \lambda \in \mathbb{R}^{d+1}_+, \qquad (2.10)$$

where $\mu_\alpha$ is the measure on $\mathbb{R}^{d+1}_+$ given by the relation (2.1).

Some basic properties of this transform are as follows. For the proofs, we refer [11, 12].

**Proposition 2.3.**   1. For all $\varphi \in L^1_\alpha(\mathbb{R}^{d+1}_+)$, the function $\mathcal{F}_{W,\alpha}(\varphi)$ is continuous on $\mathbb{R}^{d+1}_+$ and we have

$$\|\mathcal{F}_{W,\alpha}\varphi\|_{\alpha,\infty} \leq \|\varphi\|_{\alpha,1}. \qquad (2.11)$$

2. The Weinstein transform is a topological isomorphism from $\mathcal{S}_*(\mathbb{R}^{d+1}_+)$ onto itself. The inverse transform is given by

$$\mathcal{F}^{-1}_{W,\alpha}\varphi(\lambda) = \mathcal{F}_{W,\alpha}\varphi(-\lambda), \text{ for all } \lambda \in \mathbb{R}^{d+1}_+. \qquad (2.12)$$

3. Parseval formula: For all $\varphi, \phi \in \mathcal{S}_*(\mathbb{R}^{d+1}_+)$, we have

$$\int_{\mathbb{R}^{d+1}_+} \varphi(x)\overline{\phi(x)}d\mu_\alpha(x) = \int_{\mathbb{R}^{d+1}_+} \mathcal{F}_{W,\alpha}(\varphi)(x)\overline{\mathcal{F}_{W,\alpha}(\phi)(x)}d\mu_\alpha(x). \qquad (2.13)$$

4. **Plancherel formula**: For all $\varphi \in \mathcal{S}_*(\mathbb{R}^{d+1}_+)$, we have

$$\|\mathcal{F}_{W,\alpha}\varphi\|_{\alpha,2} = \|\varphi\|_{\alpha,2}. \qquad (2.14)$$

5. **Plancherel Theorem**: The Weinstein transform $\mathcal{F}_{W,\alpha}$ extends uniquely to an isometric isomorphism on $L^2_\alpha(\mathbb{R}^{d+1}_+)$.

6. **Inversion formula**: Let $\varphi \in L^1_\alpha(\mathbb{R}^{d+1}_+)$ such that $\mathcal{F}_{W,\alpha}\varphi \in L^1_\alpha(\mathbb{R}^{d+1}_+)$, then we have

$$\varphi(\lambda) = \int_{\mathbb{R}^{d+1}_+} \mathcal{F}_{W,\alpha}\varphi(x)\Lambda^d_\alpha(-\lambda, x)d\mu_\alpha(x), \ a.e. \ \lambda \in \mathbb{R}^{d+1}_+. \qquad (2.15)$$

## 2.3. The translation operator associated with the Weinstein operator

**Definition 2.4.** The translation operator $\tau^\alpha_x$, $x \in \mathbb{R}^{d+1}_+$ associated with the Weinstein operator $\Delta^d_{W,\alpha}$, is defined for a continuous function $\varphi$ on $\mathbb{R}^{d+1}_+$ which is even with respect to the last variable and for all $y \in \mathbb{R}^{d+1}_+$ by

$$\tau^\alpha_x \varphi(y) = C_\alpha \int_0^\pi \varphi\left(x' + y', \sqrt{x^2_{d+1} + y^2_{d+1} + 2x_{d+1}y_{d+1}\cos\theta}\right)(\sin\theta)^{2\alpha} d\theta,$$

with

$$C_\alpha = \frac{\Gamma(\alpha + 1)}{\sqrt{\pi}\Gamma(\alpha + 1/2)}.$$

By using the Weinstein kernel, we can also define a generalized translation, for a function $\varphi \in \mathcal{S}_*(\mathbb{R}^{d+1}_+)$ and $y \in \mathbb{R}^{d+1}_+$ the generalized translation $\tau^\alpha_x \varphi$ is defined by the following relation

$$\mathcal{F}_{W,\alpha}(\tau^\alpha_x \varphi)(y) = \Lambda^d_\alpha(x, y)\mathcal{F}_{W,\alpha}(\varphi)(y). \qquad (2.16)$$

The following proposition summarizes some properties of the Weinstein translation operator.

**Proposition 2.5.** *The translation operator $\tau_x^\alpha$, $x \in \mathbb{R}_+^{d+1}$ satisfies the following properties.*

*i). For $\varphi \in \mathbb{C}_*(\mathbb{R}^{d+1})$, we have for all $x, y \in \mathbb{R}_+^{d+1}$*

$$\tau_x^\alpha \varphi(y) = \tau_y^\alpha \varphi(x) \text{ and } \tau_0^\alpha \varphi = \varphi.$$

*ii). Let $\varphi \in L_\alpha^p(\mathbb{R}_+^{d+1})$, $1 \leq p \leq \infty$ and $x \in \mathbb{R}_+^{d+1}$. Then $\tau_x^\alpha \varphi$ belongs to $L_\alpha^p(\mathbb{R}_+^{d+1})$ and we have*

$$\|\tau_x^\alpha \varphi\|_{\alpha,p} \leq \|\varphi\|_{\alpha,p}. \tag{2.17}$$

Note that the $\mathcal{A}_\alpha(\mathbb{R}_+^{d+1})$ is contained in the intersection of $L_\alpha^1(\mathbb{R}_+^{d+1})$ and $L_\alpha^\infty(\mathbb{R}_+^{d+1})$ and hence is a subspace of $L_\alpha^2(\mathbb{R}_+^{d+1})$. For $\varphi \in \mathcal{A}_\alpha(\mathbb{R}_+^{d+1})$ we have

$$\tau_x^\alpha \varphi(y) = C_{\alpha,d} \int_{\mathbb{R}_+^{d+1}} \Lambda_\alpha^d(x,z) \Lambda_\alpha^d(-y,z) \mathcal{F}_{W,\alpha} \varphi(z) d\mu_\alpha(z). \tag{2.18}$$

By using the generalized translation, we define the generalized convolution product $\varphi *_W \psi$ of the functions $\varphi, \psi \in L_\alpha^1(\mathbb{R}_+^{d+1})$ as follows

$$\varphi *_W \psi(x) = \int_{\mathbb{R}_+^{d+1}} \tau_x^\alpha \varphi(-y) \psi(y) d\mu_\alpha(y). \tag{2.19}$$

This convolution is commutative and associative, and it satisfies the following properties.

**Proposition 2.6.** *i) For all $\varphi, \psi \in L_\alpha^1(\mathbb{R}_+^{d+1})$, (resp. $\varphi, \psi \in \mathcal{S}_*(\mathbb{R}_+^{d+1})$), then $\varphi *_W \psi \in L_\alpha^1(\mathbb{R}_+^{d+1})$, (resp. $\varphi *_W \psi \in \mathcal{S}_*(\mathbb{R}_+^{d+1})$) and we have*

$$\mathcal{F}_{W,\alpha}(\varphi *_W \psi) = \mathcal{F}_{W,\alpha}(\varphi) \mathcal{F}_{W,\alpha}(\psi). \tag{2.20}$$

*ii) Let $p, q, r \in [1, \infty]$, such that $\frac{1}{p} + \frac{1}{q} - \frac{1}{r} = 1$. Then for all $\varphi \in L_\alpha^p(\mathbb{R}_+^{d+1})$ and $\psi \in L_\alpha^q(\mathbb{R}_+^{d+1})$ the function $\varphi *_W \psi$ belongs to $L_\alpha^r(\mathbb{R}_+^{d+1})$ and we have*

$$\|\varphi *_W \psi\|_{\alpha,r} \leq \|\varphi\|_{\alpha,p} \|\psi\|_{\alpha,q}. \tag{2.21}$$

*iii) Let $\varphi, \psi \in L_\alpha^2(\mathbb{R}_+^{d+1})$. Then*

$$\varphi *_W \psi = \mathcal{F}_{W,\alpha}^{-1} \left( \mathcal{F}_{W,\alpha}(\varphi) \mathcal{F}_{W,\alpha}(\psi) \right). \tag{2.22}$$

*iv) Let $\varphi, \psi \in L_\alpha^2(\mathbb{R}_+^{d+1})$. Then $\varphi *_W \psi$ belongs to $L_\alpha^2(\mathbb{R}_+^{d+1})$ if and only if $\mathcal{F}_{W,\alpha}(\varphi) \mathcal{F}_{W,\alpha}(\psi)$ belongs to $L_\alpha^2(\mathbb{R}_+^{d+1})$ and we have*

$$\mathcal{F}_{W,\alpha}(\varphi *_W \psi) = \mathcal{F}_{W,\alpha}(\varphi) \mathcal{F}_{W,\alpha}(\psi). \tag{2.23}$$

*v) Let $\varphi, \psi \in L_\alpha^2(\mathbb{R}_+^{d+1})$. Then*

$$\|\varphi *_W \psi)\|_{\alpha,2} = \|\mathcal{F}_{W,\alpha}(\varphi) \mathcal{F}_{W,\alpha}(\psi)\|_{\alpha,2}, \tag{2.24}$$

*where both sides are finite or infinite.*

## 3. Heisenberg-Pauli-Weyl uncertainty principle

In this section we establish Heisenberg-Pauli-Weyl uncertainty principle for the operator $\mathcal{T}_{w,m,\sigma}$.

**Theorem 3.1.** *Let $m$ be a function in $L^2_\alpha(\mathbb{R}^{d+1}_+)$ satisfying the admissibility condition (1.3). Then, for $\varphi \in L^2_\alpha(\mathbb{R}^{d+1}_+)$, we have*

$$\|\varphi\|_{\alpha,2} \leq \frac{2\||y|\mathcal{F}_{W,\alpha}(\varphi)\|_{\alpha,2}}{2\alpha+d+2} \left( \int_{\mathbb{R}^{d+1}_+} \int_0^\infty |x|^2 |\mathcal{T}_{w,m,\sigma}\varphi(x)|^2_{\alpha,2} \frac{d\sigma}{\sigma} d\mu_\alpha(x) \right)^{\frac{1}{2}}. \tag{3.1}$$

*Proof.* Let $\varphi \in L^2_\alpha(\mathbb{R}^{d+1}_+)$. The inequality (3.1) holds if

$$\||y|\mathcal{F}_{W,\alpha}(\varphi)\|_{\alpha,2} = +\infty$$

or

$$\int_{\mathbb{R}^{d+1}_+} \int_0^\infty |x|^2 |\mathcal{T}_{w,m,\sigma}\varphi(x)|^2_{\alpha,2} \frac{d\sigma}{\sigma} d\mu_\alpha(x) = +\infty.$$

Let us now assume that

$$\||y|\mathcal{F}_{W,\alpha}(\varphi)\|_{\alpha,2} + \int_{\mathbb{R}^{d+1}_+} \int_0^\infty |x|^2 |\mathcal{T}_{w,m,\sigma}\varphi(x)|^2_{\alpha,2} \frac{d\sigma}{\sigma} d\mu_\alpha(x) < +\infty.$$

Inequality (1.2) leads to

$$\int_{\mathbb{R}^{d+1}_+} |\mathcal{T}_{w,m,\sigma}\varphi(x)|^2_{\alpha,2} d\mu_\alpha(x)$$

$$< \left( \int_{\mathbb{R}^{d+1}_+} |x|^2 |\mathcal{T}_{w,m,\sigma}\varphi(x)|^2_{\alpha,2} d\mu_\alpha(x) \right)^{\frac{1}{2}}$$

$$\times \left( \int_{\mathbb{R}^{d+1}_+} |y|^2 |\mathcal{F}_{W,\alpha}(\mathcal{T}_{w,m,\sigma}\varphi(.))(y)|^2_{\alpha,2} d\mu_\alpha(y) \right)^{\frac{1}{2}}.$$

Integrating with respect to $d\sigma/\sigma$, we get

$$\|\mathcal{T}_{w,m,\sigma}\varphi)\|^2_{\alpha,2} < \int_0^\infty \left( \int_{\mathbb{R}^{d+1}_+} |x|^2 |\mathcal{T}_{w,m,\sigma}\varphi(x)|^2_{\alpha,2} d\mu_\alpha(x) \right)^{\frac{1}{2}}$$

$$\times \left( \int_{\mathbb{R}^{d+1}_+} |y|^2 |\mathcal{F}_{W,\alpha}(\mathcal{T}_{w,m,\sigma}\varphi(.))(y)|^2_{\alpha,2} d\mu_\alpha(y) \right)^{\frac{1}{2}} \frac{d\sigma}{\sigma}.$$

From [14, Theorem 2.3] and Schwarz's inequality, we obtain

$$\|\varphi\|_{\alpha,2}^2 < \left(\int_0^\infty \int_{\mathbb{R}_+^{d+1}} |x|^2 |\mathcal{T}_{w,m,\sigma}\varphi(x)|_{\alpha,2}^2 d\mu_\alpha(x) \frac{d\sigma}{\sigma}\right)^{\frac{1}{2}}$$

$$\times \left(\int_0^\infty \int_{\mathbb{R}_+^{d+1}} |y|^2 |\mathcal{F}_{W,\alpha}(\mathcal{T}_{w,m,\sigma}\varphi(.))(y)|_{\alpha,2}^2 d\mu_\alpha(y) \frac{d\sigma}{\sigma}\right)^{\frac{1}{2}}.$$

From (1.1), Fubini-Tonnelli's theorem and the admissibility condition (1.3), we have

$$\int_0^\infty \int_{\mathbb{R}_+^{d+1}} |y|^2 |\mathcal{F}_{W,\alpha}(\mathcal{T}_{w,m,\sigma}\varphi(.))(y)|_{\alpha,2}^2 d\mu_\alpha(y) \frac{d\sigma}{\sigma}$$
$$= \int_0^\infty \int_{\mathbb{R}_+^{d+1}} |y|^2 |m_\sigma(y)|^2 |\mathcal{F}_{W,\alpha}(\varphi)(y)|_{\alpha,2}^2 d\mu_\alpha(y) \frac{d\sigma}{\sigma}$$
$$= \int_{\mathbb{R}_+^{d+1}} |y|^2 |\mathcal{F}_{W,\alpha}(\varphi)(y)|_{\alpha,2}^2 d\mu_\alpha(y).$$

This gives the result and completes the proof of the theorem. □

**Theorem 3.2.** *Let $m$ be a function in $L_\alpha^2(\mathbb{R}_+^{d+1})$ satisfying the admissibility condition (1.3) and $\beta, \delta \in [1, \infty)$. Let $\varepsilon \in \mathbb{R}$, such that $\beta\varepsilon = (1-\varepsilon)\delta$ then, for $\varphi \in L_\alpha^2(\mathbb{R}_+^{d+1})$, we have*

$$\|\varphi\|_{\alpha,2} \leq \left(\frac{2}{2\alpha+d+2}\right)^{\beta\varepsilon} \left\||x|^\beta \mathcal{T}_{w,m,\sigma}\varphi\right\|_{\alpha,2}^\varepsilon \left\||y|^\delta \mathcal{F}_{W,\alpha}(\varphi)\right\|_{\alpha,2}^{1-\varepsilon}. \qquad (3.2)$$

*Proof.* Let $\varphi \in L_\alpha^2(\mathbb{R}_+^{d+1})$. The inequality (3.1) holds if

$$\left\||x|^\beta \mathcal{T}_{w,m,\sigma}\varphi\right\|_{\alpha,2}^\varepsilon = +\infty \quad \text{or} \quad \left\||y|^\delta \mathcal{F}_{W,\alpha}(\varphi)\right\|_{\alpha,2}^{1-\varepsilon} = +\infty.$$

Let us now assume that $\varphi \in L_\alpha^2(\mathbb{R}_+^{d+1})$ with $\varphi \neq 0$ such that

$$\left\||x|^\beta \mathcal{T}_{w,m,\sigma}\varphi\right\|_{\alpha,2}^\varepsilon + \left\||y|^\delta \mathcal{F}_{W,\alpha}(\varphi)\right\|_{\alpha,2}^{1-\varepsilon} < +\infty,$$

therefore, for all $\delta > 1$, we have

$$\left\||x|^\beta \mathcal{T}_{w,m,\sigma}\varphi\right\|_{\alpha,2}^{\frac{1}{\beta}} \|\mathcal{T}_{w,m,\sigma}\varphi\|_{\alpha,2}^{\frac{1}{\beta'}} = \left\||x|^2 |\mathcal{T}_{w,m,\sigma}\varphi|^{\frac{2}{\beta}}\right\|_{\alpha,\beta}^{\frac{1}{2}} \left\||\mathcal{T}_{w,m,\sigma}\varphi|^{\frac{2}{\beta'}}\right\|_{\alpha,\beta'}^{\frac{1}{2}},$$

with $\beta' = \frac{\beta}{\beta-1}$.

Applying the Hölder's inequality, we get

$$\||x|\mathcal{T}_{w,m,\sigma}\varphi\|_{\alpha,2} \leq \left\||x|^\beta \mathcal{T}_{w,m,\sigma}\varphi\right\|_{\alpha,2}^{\frac{1}{\beta}} \|\mathcal{T}_{w,m,\sigma}\varphi\|_{\alpha,2}^{\frac{1}{\beta'}}.$$

According to [14, Theorem 2.3], we have for all $\beta \geq 1$

$$\||x|\mathcal{T}_{w,m,\sigma}\varphi\|_{\alpha,2} \leq \left\||x|^\beta \mathcal{T}_{w,m,\sigma}\varphi\right\|_{\alpha,2}^{\frac{1}{\beta}} \|\varphi\|_{\alpha,2}^{\frac{1}{\beta'}}, \qquad (3.3)$$

with equality if $\beta = 1$. In the same manner, for all $\delta \geq$ and using Plancherel formula (2.14), we get

$$\| |y| \mathcal{F}_{W,\alpha}(\varphi) \|_{\alpha,2} \leq \| |y|^\delta \mathcal{F}_{W,\alpha}(\varphi) \|_{\alpha,2}^{\frac{1}{\delta}} \| \varphi \|_{\alpha,2}^{\frac{1}{\delta'}}, \tag{3.4}$$

with equality if $\delta = 1$. By using the fact that $\beta\varepsilon = (1-\varepsilon)\delta$ and according to inequalities (3.3) and (3.4), we have

$$\left( \frac{\| |x| \mathcal{T}_{w,m,\sigma}\varphi \|_{\alpha,2} \| |y| \mathcal{F}_{W,\alpha}(\varphi) \|_{\alpha,2}}{\| \varphi \|_{\alpha,2}^{\frac{1}{\beta'}+\frac{1}{\delta'}}} \right)^{\beta\delta}$$
$$\leq \| |x|^\beta \mathcal{T}_{w,m,\sigma}\varphi \|_{\alpha,2}^\varepsilon \| |y|^\delta \mathcal{F}_{W,\alpha}(\varphi) \|_{\alpha,2}^{1-\varepsilon},$$

with equality if $\beta = \delta = 1$. Next by Theorem 3.1, we obtain

$$\| \varphi \|_{\alpha,2} \leq \left( \frac{2}{2\alpha+d+2} \right)^{\beta\varepsilon} \| |x|^\beta \mathcal{T}_{w,m,\sigma}\varphi \|_{\alpha,2}^\varepsilon \| |y|^\delta \mathcal{F}_{W,\alpha}(\varphi) \|_{\alpha,2}^{1-\varepsilon},$$

which completes the proof of the theorem. $\square$

## 4. Donoho-Stark's uncertainty principle

**Definition 4.1.** (i) Let $\Omega$ be a measurable subset of $\mathbb{R}_+^{d+1}$, we say that the function $\varphi \in L_\alpha^2(\mathbb{R}_+^{d+1})$ is $\epsilon$-concentrated on $\Omega$, if

$$\| \varphi - \chi_\Omega \varphi \|_{\alpha,2} \leq \epsilon \| \varphi \|_{\alpha,2}, \tag{4.1}$$

where $\chi_\Omega$ is the indicator function of the set $\Omega$.

(ii) Let $\Sigma$ be a measurable subset of $(0,\infty) \times \mathbb{R}_+^{d+1}$ and let $\varphi \in L_\alpha^2(\mathbb{R}_+^{d+1})$. We say that $\mathcal{T}_{w,m,\sigma}\varphi$ is $\nu$-concentrated on $\Sigma$, if

$$\| \mathcal{T}_{w,m,\sigma}\varphi - \chi_\Sigma \mathcal{T}_{w,m,\sigma}\varphi \|_{2,\alpha} \leq \nu \| \mathcal{T}_{w,m,\sigma} \|_{2,\alpha}, \tag{4.2}$$

where $\chi_\Sigma$ is the indicator function of the set $\Sigma$.

We need the following Lemma for the proof of Donoho-Stark's uncertainty principle.

**Lemma 4.2.** Let $m, \varphi \in L_\alpha^1(\mathbb{R}_+^{d+1}) \cap L_\alpha^2(\mathbb{R}_+^{d+1})$. Then the operators $\mathcal{T}_{w,m,\sigma}$ satisfy the following integral representation.

$$\mathcal{T}_{w,m,\sigma} = \frac{1}{\sigma^{2\alpha+d+2}} \int_{\mathbb{R}_+^{d+1}} \Psi_\alpha(x,y) \varphi(y) d\mu_\alpha(y), \quad (\sigma, x) \in (0,\infty) \times \mathbb{R}_+^{d+1},$$

where

$$\Psi_\alpha(x,y) = \int_{\mathbb{R}_+^{d+1}} m_\sigma(z) \Lambda_\alpha^d(\lambda, x) \Lambda_\alpha^d(\lambda, -y) d\mu_\alpha(z).$$

*Proof.* The result follows from the definition of the Weinstein $L^2$-Multiplier operators (1.1) and the inversion formula of the Weinstein transform (2.12) using Fubini-Tonnelli's theorem. $\square$

**Theorem 4.3.** *Let $\varphi$ be a function in $L^2_\alpha(\mathbb{R}^{d+1}_+)$ and $m \in L^1_\alpha(\mathbb{R}^{d+1}_+) \cap L^2_\alpha(\mathbb{R}^{d+1}_+)$ satisfying the admissibility condition (1.3). If $\varphi$ is $\epsilon$-concentrated on $\Omega$ and $\mathcal{T}_{w,m,\sigma}\varphi$ is $\nu$-concentrated on $\Sigma$, then*

$$\|m\|_{\alpha,1} (\mu_\alpha(\Omega))^{\frac{1}{2}} \left( \int\int_\Sigma \frac{1}{\sigma^{2(2\alpha+d+2)}} d\Theta_\alpha(\sigma,x) \right)^{\frac{1}{2}} \geq 1 - (\epsilon + \nu),$$

*where $\Theta_\alpha$ is the measure on $(0,\infty) \times \mathbb{R}^{d+1}_+$ given by $d\Theta_\alpha(\sigma,x) := (d\sigma/\sigma)d_\alpha\mu(x)$.*

*Proof.* Let $\varphi$ be a function in $L^2_\alpha(\mathbb{R}^{d+1}_+)$. Assume that $0 < \mu_\alpha(\Omega) < \infty$ and

$$\int\int_\Sigma \frac{1}{\sigma^{2(2\alpha+d+2)}} d\Theta_\alpha(\sigma,x) < \infty.$$

According to [14, Theorem 2.3] and inequalities (4.1)-(4.2), we get

$$\begin{aligned}
\|\mathcal{T}_{w,m,\sigma}\varphi - \chi_\Sigma \mathcal{T}_{w,m,\sigma}(\chi_\Omega\varphi)\|_{2,\alpha} &\leq \|\mathcal{T}_{w,m,\sigma}\varphi - \chi_\Sigma \mathcal{T}_{w,m,\sigma}\varphi\|_{2,\alpha} \\
&\quad + \|\mathcal{T}_{w,m,\sigma}\varphi - \chi_\Sigma \mathcal{T}_{w,m,\sigma}(\varphi - \chi_\Omega\varphi)\|_{2,\alpha} \\
&\leq \|\mathcal{T}_{w,m,\sigma}\varphi - \chi_\Sigma \mathcal{T}_{w,m,\sigma}(\varphi - \chi_\Omega\varphi)\|_{2,\alpha} \\
&\quad + \nu\|\mathcal{T}_{w,m,\sigma}\varphi\|_{2,\alpha} \\
&\leq (\epsilon + \nu)\|\varphi\|_{2,\alpha}.
\end{aligned}$$

By triangle inequality it follows that

$$\begin{aligned}
\|\mathcal{T}_{w,m,\sigma}\varphi\|_{2,\alpha} &\leq \|\mathcal{T}_{w,m,\sigma}\varphi - \chi_\Sigma \mathcal{T}_{w,m,\sigma}(\chi_\Omega\varphi)\|_{2,\alpha} + \|\chi_\Sigma \mathcal{T}_{w,m,\sigma}(\chi_\Omega\varphi)\|_{2,\alpha} \\
&\leq (\epsilon + \nu)\|\varphi\|_{2,\alpha} + \|\chi_\Sigma \mathcal{T}_{w,m,\sigma}(\chi_\Omega\varphi)\|_{2,\alpha}. \quad (4.3)
\end{aligned}$$

On the other hand, we have

$$\|\chi_\Sigma \mathcal{T}_{w,m,\sigma}(\chi_\Omega\varphi)\|_{2,\alpha} = \left( \int\int_\Sigma |\mathcal{T}_{w,m,\sigma}(\chi_\Omega\varphi)(x)|^2 d\Theta_\alpha(\sigma,x) \right)^{\frac{1}{2}}$$

and moreover $m, \chi_\Omega\varphi \in L^1_\alpha(\mathbb{R}^{d+1}_+) \cap L^2_\alpha(\mathbb{R}^{d+1}_+)$, then by Lemma 4.2, we obtain

$$|\mathcal{T}_{w,m,\sigma}(\chi_\Omega\varphi)(x)| \leq \frac{1}{\sigma^{2\alpha+d+2}} \|m\|_{1,\alpha} \|\varphi\|_{2,\alpha} (\mu_\alpha(\Omega))^{\frac{1}{2}}.$$

Therefore, thus

$$\begin{aligned}
\|\chi_\Sigma \mathcal{T}_{w,m,\sigma}(\chi_\Omega\varphi)\|_{2,\alpha} &\leq \|m\|_{1,\alpha}\|\varphi\|_{2,\alpha}(\mu_\alpha(\Omega))^{\frac{1}{2}} \\
&\quad \times \left( \int\int_\Sigma \frac{1}{\sigma^{2(2\alpha+d+2)}} d\Theta_\alpha(\sigma,x) \right)^{\frac{1}{2}}.
\end{aligned}$$

Hence, according to last inequality and (4.3)

$$\begin{aligned}
\|\mathcal{T}_{w,m,\sigma}(\varphi)\|_{2,\alpha} &\leq \|m\|_{1,\alpha}\|\varphi\|_{2,\alpha}(\mu_\alpha(\Omega))^{\frac{1}{2}} \\
&\quad \times \left( \int\int_\Sigma \frac{1}{\sigma^{2(2\alpha+d+2)}} d\Theta_\alpha(\sigma,x) \right)^{\frac{1}{2}} + (\epsilon + \nu)\|\varphi\|_{2,\alpha}.
\end{aligned}$$

Applying Plancherel formula [14, Theorem 2.3], we obtain

$$\|m\|_{\alpha,1} (\mu_\alpha(\Omega))^{\frac{1}{2}} \left( \int\int_\Sigma \frac{1}{\sigma^{2(2\alpha+d+2)}} d\Theta_\alpha(\sigma,x) \right)^{\frac{1}{2}} \geq 1 - (\epsilon + \nu),$$

which completes the proof of the theorem. □

**Corollary 4.4.** *If $\Sigma = \{(\sigma, x) \in (0, \infty) \times \mathbb{R}_+^{d+1} : \sigma \geq \varrho\}$ for some $\varrho > 0$, one assumes that*

$$\rho = \max \left\{ 1/\sigma : (\sigma, x) \in \Sigma \text{ for some } x \in \mathbb{R}_+^{d+1} \right\}.$$

*Then by the previous Theorem, we deduce that*

$$\rho^{2\alpha+d+2} \|m\|_{\alpha,1} \left(\mu_\alpha(\Omega)\right)^{\frac{1}{2}} \left(\Theta_\alpha(\Sigma)\right)^{\frac{1}{2}} \geq 1 - (\epsilon + \nu).$$

Ahmed Saoudi

Northern Border University, College of Science, Arar, P.O. Box 1631, Saudi Arabia.
Université de Tunis El Manar, Faculté des sciences de Tunis,
LR13ES06 Laboratoire de Fonctions spéciales, analyse harmonique et analogues,
2092 Tunis, Tunisia
e-mail: `ahmed.saoudi@ipeim.rnu.tn`